\documentclass{amsart}
\usepackage{amssymb,latexsym}
\usepackage{isolatin1}
\numberwithin{equation}{section}

\newtheorem{Lem}{Lemma}[section]

\newtheorem{Cor}[Lem]{Corollary}
\newtheorem{Thm}[Lem]{Theorem}

\theoremstyle{definition}
\newtheorem{Def}[Lem]{Definition}
\theoremstyle{remark}
\newtheorem{Rem}[Lem]{Remark}

\renewcommand\o{\otimes}

\DeclareMathOperator\ev{\operatorname{ev}}
\DeclareMathOperator\db{\operatorname{db}}
\newcommand\op{{\operatorname{op}}}
\newcommand\cop{{\operatorname{cop}}}
\newcommand\bop{{\operatorname{bop}}}

\newcommand\HMod[4]{{^{#1}_{#3}\mathcal M^{#2}_{#4}}}
\newcommand\LMod[1]{{_{#1}\mathcal M}}

\newcommand\BiMod[1]{{_{#1}\mathcal M_{#1}}}

\newcommand{\co}[1]{\mathrel{\mathop{\Box}_{#1}}}

\newcommand\sw[1]{{}_{(#1)}}

\newcommand\so[1]{^{(#1)}}
\newcommand\som[1]{^{(-#1)}}

\newcommand\inv{^{-1}}

\renewcommand\epsilon\varepsilon

\newcommand\role{r\^ole}
\makeatletter
\def\namelabel#1#2{\@bsphack
  \protected@write\@auxout{}%
         {\string\newlabel{#1.nme}{{#2}{#2}}}%
  \@esphack}
\def\nmlabel#1#2{\label{#2}\namelabel{#2}{#1}}
\newcommand\nmref[1]{\ref{#1.nme}\ \ref{#1}}
\makeatother

\begin{document}
\title{A quasi-Hopf algebra freeness theorem}
\author{Peter Schauenburg}
\address{Mathematisches Institut der Universit\"at M\"unchen, 
Theresienstr.~39, 80333~M\"unchen, Germany}
\email{schauen@rz.mathematik.uni-muenchen.de}
\begin{abstract}
  We prove the quasi-Hopf algebra version of the 
  Nichols-Zoeller theorem: A finite-dimensional quasi-Hopf algebra
  is free over any quasi-Hopf subalgebra.
\end{abstract}
\maketitle
\section{Introduction}
Kaplansky's first conjecture \cite{Kap:B} asked whether a 
Hopf algebra $H$ is always free over any Hopf subalgebra $K$.
While this is false in the infinite-dimensional case, as shown
by Oberst and Schneider \cite{ObeSch:UFGEI}, it is true when 
$H$ is finite-dimensional. This is the content of 
the celebrated Nichols-Zoeller theorem \cite{NicZoe:HAFT}.
More generally, Nichols and Zoeller prove that Hopf modules in 
the category $\HMod{}H{}K$ are always free as $K$-modules, provided
again that $H$ is a finite-dimensional Hopf algebra and $K$ a 
Hopf subalgebra. The Nichols-Zoeller theorem is an indispensable
tool in the study of finite-dimensional Hopf algebras.

In the present paper, we will prove the essentially the same results
for quasi-Hopf algebras in place of Hopf algebras.
Quasi-Hopf algebras
were introduced by Drinfeld \cite{Dri:QHA}. By definition, 
a quasibialgebra
$H$ is an algebra and a non-coassociative coalgebra, whose
lack of coassociativity is controlled by an invertible element
$\phi\in H\o H\o H$, the associator; more precisely, comultiplication is coassociative
up to conjugation by $\phi$. The ``meaning'' of the definition is
already explained in Drinfeld's paper: The axioms are such that
the category $\LMod H$ of modules over $H$ is a monoidal category
with respect to the diagonal module structure on the tensor product
of $H$-modules, taken over the base field $k$. The rather complicated
axioms for a quasiantipode are designed so that one can define
a dual $H$-module for any finite-dimensional $H$-module.

The question for a Nichols-Zoeller theorem for quasi-Hopf algebras
was brought to the author's
attention by Robert Fischer.
We shall prove that a quasi-Hopf algebra $H$ is a free 
module over any
quasi-Hopf subalgebra $K$. This is by definition a subalgebra
and (non-coassociative) subcoalgebra such that $K\o K\o K$ contains
the associator $\phi$. Interesting examples of quasi-Hopf algebra
inclusions in this sense have recently been studied by 
Masuoka \cite{Mas:HAEC,Mas:CCBEAMPB}. We note that in Masuoka's 
examples, which are generalizations of bicrossproduct constructions,
$H$ is a free $K$-module by construction. The same is true for 
the embedding of a finite-dimensional quasi-Hopf algebra in its
Drinfeld double, cf.\ \cite{HauNil:DCPDQQG,HauNil:DQQG}.

We will also prove a quasi-Hopf version of the freeness result
of Nichols and Zoeller on Hopf modules.
However, we will not
use this result to prove freeness of $H$ over $K$. In fact the 
freeness of $H$ over $K$ will be proved in the more general situation
that $K$ is a subalgebra of $H$, which has some structure of
quasi-Hopf algebra with 
its own associator that needs not coincide with that of $H$.
Natural examples of such a situation come immediately
to mind: If $K$ and $F$ are
quasi-Hopf algebras, then $K\o F:=H$ is a quasi-Hopf algebra, 
containing
$K$ as a subalgebra and subcoalgebra, but with its associator
not contained in $K\o K\o K$, unless the associator of $F$
is trivial. 
(Incidentally, $H$ is of course
a free $K$-module in this situation.)

Our proofs will
follow the original proofs of Nichols and Zoeller, or more 
precisely the account in Montgomery's book \cite{Mon:HAAR}, 
quite closely.
A lemma on modules over Frobenius algebras, in particular, 
will retain its key \role\ in the proof. The main obstacle for the
generalization is the constant use of Hopf module categories
like $\HMod{}H{}K$, and of structure theorems for Hopf modules.
Note that the said Hopf module category is not defined when
$H$ is a quasi-Hopf algebra, since then $H$ need not be a coalgebra.

These difficulties can be overcome by considering suitable 
Hopf module categories that have been defined by Hausser and Nill
\cite{HauNil:ITQHA}. The point is that while a quasibialgebra is not
a coalgebra, one can still define Hopf bimodule categories 
$\HMod{}HHH$. This is surprising at first sight, but the definition
by Hausser and Nill turns out to be simply the definition of an 
$H$-comodule --- but within the monoidal category $\BiMod H$ of 
$H$-bimodules, in which $H$ {\em is} a coalgebra after all.
This categorical viewpoint on Hopf modules was stressed and used 
in \cite{Sch:HMDQHA}.
Hausser and Nill also proved a structure theorem for
Hopf modules in $\HMod{}HHH$, generalizing the well-known structure
theorem for Hopf modules in $\HMod{}H{}H$ with an ordinary Hopf 
algebra $H$.
\section{Preliminaries}
Throughout the paper, we work over some base field $k$. Tensor products,
algebras and the like are over $k$. 
\subsection{Quasi-Hopf algebras}
We recall the definition and some basic facts concerning quasibialgebras,
quasi-Hopf algebras, and their relation to monoidal categories. 
General references are
Drinfeld's original paper \cite{Dri:QHA}, and Kassel's book
\cite{Kas:QG}.

A quasibialgebra $H=(H,\Delta,\epsilon,\phi)$ consists of an algebra $H$,
algebra maps $\Delta\colon H\rightarrow H\o H$ and $\epsilon\colon H\rightarrow k$,
and an invertible element $\phi\in H^{\o 3}$, the associator, such that
\begin{gather}
(\epsilon\o H)\Delta(h)=h=(H\o\epsilon)\Delta(h),\\
(H\o\Delta)\Delta(h)\cdot\phi=\phi\cdot(\Delta\o H)\Delta(h)\label{quasicoass},\\
(H\o H\o\Delta)(\phi)\cdot(\Delta\o H\o H)(\phi)=(1\o\phi)\cdot(H\o\Delta\o H)(\phi)\cdot(\phi\o 1)\label{cocycle},\\
(H\o\epsilon\o H)(\phi)=1
\end{gather}
hold for all $h\in H$.
We will write $\Delta(h)=:h\sw 1\o h\sw 2$, 
$\phi=\phi\so 1\o\phi\so 2\o\phi\so 3$, and
$\phi\inv=\phi\som 1\o\phi\som 2\o\phi\som 3$. 

If $H$ is a quasibialgebra, then the category $\LMod H$ of
left $H$-modules is a monoidal category in the following way:
For $V,W\in\LMod H$ the tensor product $V\o W$ is an $H$-module
by $h(v\o w)=h\sw 1v\o h\sw 2w$. The base
ring $k$ is an $H$-module via $\epsilon$. The canonical
morphisms $V\cong V\o k\cong k\o V$ are $H$-linear for $V\in\LMod H$.
The map
$$\Phi\colon (U\o V)\o W\ni u\o v\o w\mapsto\phi\so 1u\o\phi\so 2v\o\phi\so 3w\in U\o (V\o W)$$
for $H$-modules $U,V,W$ is $H$-linear as a consequence of 
\eqref{quasicoass}, and satisfies 
Mac Lane's pentagon axiom for a monoidal category as a 
consequence of \eqref{cocycle}.

We denote $H$ with the opposite multiplication, coopposite comultiplication,
or both, by $H^\op$, $H^\cop$, and $H^\bop$, respectively.
If $(H,\phi)$ is a quasibialgebra, so is $(H^\op,\phi\inv)$,
$(H^\cop,\phi\so 3\o\phi\so 2\o\phi\so 1)$, and hence
also $(H^\bop,\phi\som 3\o\phi\som 2\o\phi\som 1)$.
If $(B,\psi)$ is another quasibialgebra, then 
$H\o B$ is a quasibialgebra with associator
$\phi\so 1\o\psi\so 1\o\phi\so 2\o\psi\so 2\o\phi\so 3\o\psi\so 3$.

A quasiantipode $(S,\alpha,\beta)$ for a quasibialgebra $H$ consists
of an algebra antiautomorphism $S$ of $H$, and elements $\alpha,\beta\in H$,
such that 
\begin{align*}
  S(h\sw 1)\alpha h\sw 2&=\epsilon(h)\alpha, & h\sw 1\beta S(h\sw 2)&=\epsilon(h)\beta,\\
  \phi\so 1\beta S(\phi\so 2)\alpha\phi\so 3&=1,&S(\phi\som 1)\alpha\phi\som 2\beta\phi\som 3&=1
\end{align*}
hold in $H$, for $h\in H$. A quasi-Hopf algebra is a quasibialgebra
with a quasi-antipode. If $H$ is a quasi-Hopf algebra, then so are
$H^\op$, $H^\cop$, and $H^\bop$; for later use we note that the 
underlying map of the antipode of $H^\op$ is $S\inv$. 
Let $H$ be a quasi-Hopf algebra.
If $V$ is a finite-dimensional
left $H$-module, then $V^*$, the dual vector space, becomes an 
left $H$-module with the action induced by the transpose of the 
original action of $H$ via the antipode, that is, by setting
$\langle h\cdot\varphi,v\rangle=\langle \varphi,S(h)v\rangle$,
for $h\in H$, $\varphi\in V^*$, and $v\in V$, 
where $\langle,\rangle$ denotes evaluation. 
Together with suitable
evaluation and coevaluation maps
$$\ev\colon V^*\o V\rightarrow k,\qquad\db\colon k\rightarrow V\o V^*$$
the object $V^*$ is a dual object to $V$ in the monoidal category
of $H$-modules.
\subsection{Hopf modules}
Let $H$ be a quasibialgebra. Then it is impossible to define a 
comodule over $H$, since $H$ is not a coassociative coalgebra.
It is also impossible to define a Hopf module category
$\HMod{}H{}H$ with one module and one comodule structure.

However, there is a definition of Hopf modules in $\HMod{}HHH$, which
can be viewed as a special case of the theory of (co)algebras and
(co)modules within monoidal categories. For a development of this
theory, in which many of the elementary facts and constructions
of ordinary ring theory can be carried over to quite arbitrary
abstract monoidal categories replacing the category of abelian 
groups or $k$-modules, we refer the reader to
\cite{Par:NARMTI,Par:NARMTII}. 
We have applied this point of view on Hopf modules 
over quasi-Hopf algebras in \cite{Sch:HMDQHA};
it allows us to get away 
without lengthy calculations involving
associators and the explicit use of the quasibialgebra axioms
in many cases.

The axioms of a quasibialgebra imply that $H$ is a 
coassociative coalgebra within the monoidal category 
$\BiMod H$ of $H$-$H$-bimodules: The associator morphism
in $\BiMod H\cong\LMod{H^\op\o H}$ is given by
$$\Phi\colon (U\o V)\o W\ni u\o v\o w\mapsto\phi\so 1u\phi\som 1\o \phi\so 2v\phi\som 2\o\phi\so 3w\phi\som 3\in U\o(V\o W),$$
so that coassociativity of $H$ as a coalgebra in $\BiMod H$ is in 
fact the same as \eqref{quasicoass}.

Now one can define Hopf modules over $H$, as first introduced by
Hausser and Nill \cite{HauNil:ITQHA}, as follows:
\begin{Def}
  Let $H$ be a quasibialgebra. A Hopf module $M\in\HMod{}HHH$ is
  a right $H$-comodule within the monoidal category $\BiMod H$
  of $H$-bimodules. Explicitly, $M$ is an $H$-$H$-bimodule equipped
  with a map $\delta\colon M\rightarrow M\o H$, $\delta(m)=m\sw 0\o m\sw 1$, 
  satisfying $m\sw 0\epsilon(m\sw 1)=m$ and
  $$(M\o\Delta)\delta(m)\cdot\phi=\phi\cdot(\delta\o H)\delta(m)$$
  for all $m\in M$, where the multiplication by $\phi$ now takes 
  place in the $H\o H\o H$-bimodule $M\o H\o H$.

  More generally, let $K\subset H$ be a quasi-Hopf subalgebra. Then
  $H$ is also a coalgebra in $\BiMod K$, and we define the 
  Hopf module category $\HMod{}HKK$ to be the category of $H$-comodules
  within $\BiMod K$.

  We note that Hopf module categories like $\HMod H{}KK,\HMod HHKK$
  can be defined in the same way.
\end{Def}

The interpretation of Hopf modules as comodules within some 
monoidal category gives an immediate supply of constructions
of and with such objects.  
We refer the reader to \cite{Sch:HMDQHA} for 
examples rather than giving explicit details on the following
uses of the general principle:
\begin{Rem}\nmlabel{Remark}{remzeug}
  Let $H$ be a quasibialgebra and $K$ a quasi-Hopf subalgebra. 
  \begin{enumerate}
    \item Whenever $P\in\BiMod K$, and $M\in\HMod{}HKK$, then 
      $ {_\cdot P_\cdot}\o{_\cdot M_\cdot^\cdot}\in\HMod{}HKK$.
      Here and below the dots are supposed to indicate that
      $M\o P$ is endowed with the diagonal left and right 
      $K$-module structures, and the right $H$-comodule structure
      induced by that of $P$, but taken within the monoidal category
      $\BiMod K$ with its nontrivial associator isomorphism.
      Explicitly, this means that the comodule structure is given by
      \begin{align*}
        P\o M&\rightarrow (P\o M)\o H,\\
        p\o m&\mapsto \phi\som 1 p\phi\so 1\o \phi\som 2m\sw 0\phi\so 2\o\phi\som 3m\sw 1\phi\so 3.
      \end{align*}
    \item\label{cofree} As a special case, for any $P\in\BiMod K$ we have
      ${_\cdot P_\cdot}\o {_\cdot H_\cdot^\cdot}\in\HMod{}HKK.$
      This is the cofree right $H$-comodule generated by $P$ within
      the monoidal category $\BiMod K$. 
    \item\label{speziell} Again as a special case, one
      can consider ${_\cdot V}\o{_\cdot H_\cdot^\cdot}\in\HMod{}HKK$
      for any $V\in\LMod K$.
    \item\label{cotens} As in any monoidal category, there is a notion of 
      cotensor product of a right and a left $H$-comodule within 
      $\BiMod K$, that is, a cotensor product 
      $M\co HN\in\BiMod K$ for $M\in\HMod{}HKK$ and $N\in\HMod H{}KK$.
      About this construction, we shall only need to know that
      $$M\co H({_\cdot ^\cdot H_\cdot}\o{_\cdot P_\cdot})\cong M\o P$$
      as $K$-$K$-bimodules, for any $M\in\HMod{}HKK$ and $P\in\BiMod K$.
      Here, ${_\cdot ^\cdot H_\cdot}\o{_\cdot P_\cdot}\in\HMod H{}KK$
      is a left-right switched version of the construction in 
      \eqref{cofree}.
  \end{enumerate}
\end{Rem}
The structure theorem for Hopf modules over quasi-Hopf algebras 
proved by Hausser and Nill \cite[Prop.3.11.]{HauNil:ITQHA} 
says that, for
a quasi-Hopf algebra $H$, the Hopf modules in $\HMod{}HHH$ constructed
in part \eqref{speziell} of \nmref{remzeug} are the only ones there
are:
\begin{Thm}\nmlabel{Theorem}{strukthm}
  Let $H$ be a quasi-Hopf algebra. Then the functor
  $$\LMod H\ni V\mapsto {_\cdot V}\o{_\cdot H_\cdot^\cdot}\in\HMod{}HHH$$
  is a category equivalence. In particular, every Hopf module
  in $\HMod{}HHH$ is a free right $H$-module.
\end{Thm}
We take the liberty to refer to \cite{Sch:HMDQHA} for another proof.

The main application of Hopf module theory in \cite{HauNil:ITQHA} is
the development of an integral theory for quasi-Hopf algebras.
In particular, Hausser and Nill 
 \cite[Thm.4.3]{HauNil:ITQHA} show the
 existence of
a cointegral $\lambda\in H^*$ for any finite-dimensional 
quasi-Hopf algebra $H$, which is nondegenerate and hence makes 
$H$ a Frobenius algebra. Since the fact that $H$ is Frobenius is
of key importance for the Nichols-Zoeller result, we shall sketch
a short
proof of this fact for completeness:
\begin{Thm}
Let $H$ be a finite-dimensional quasi-Hopf algebra. Then there exists
a non-degenerate linear form $\lambda\in H^*$, so that $H$ is a 
Frobenius algebra.
\end{Thm}
(Nondegeneracy means that $H\times H\ni (g,h)\mapsto \lambda(gh)\in k$
is a nondegenerate bilinear form.)
\begin{proof}
Since $H$ is a quasi-Hopf algebra, so is $H^\op$, and hence
$H^\op\o H$. 
This implies that for each $V\in\BiMod H\cong\LMod{H^\op\o H}$
the dual space $V^*$ becomes a dual object for $V$ within the monoidal
category $\BiMod H$, with module structures induced by the antipode
of $H^\op\o H$. More precisely, the $H$-$H$-bimodule structure
on $V^*$ is given by
$$\langle g\cdot \varphi\cdot h,v\rangle=\langle\varphi,S(g)vS\inv(h)\rangle$$
for $v\in V$, $\varphi\in V^*$, and $g,h\in H$, where 
$\langle,\rangle$
denotes evaluation. 
It is easy to check (and a very special case of
\cite[Cor.3.7.]{Par:NARMTI}) that when an object $V$ of a 
monoidal category $\mathcal C$ is a left comodule over some coalgebra 
$C\in\mathcal C$,
and has a dual $V^*$, then $V^*$ is a right $C$-comodule.
For example, since $H$ is a left $H$-comodule within the 
monoidal category $\BiMod H$, its dual $H^*$ is a right $H$-comodule
(whose structure map we will not need to know), that is, an object
of $\HMod{}HHH$. 
It follows from \nmref{strukthm} that there is a (one-dimensional) left
$H$-module $V$ and an isomorphism
$$\psi\colon {_\cdot V}\o{_\cdot H^\cdot_\cdot}\rightarrow H^*$$
of Hopf modules in $\HMod{}HHH$. Since $\psi$ is a right $H$-module
map, it has the form $\psi(v\o h)=\psi_0(v)h$ for some 
$\psi_0\colon V\rightarrow H^*$. Let $v\in V\setminus\{0\}$, and set
$\lambda:=\psi_0(v)$. Then it follows that
$$H\ni h\mapsto (\lambda\cdot h\colon g\mapsto \lambda(S\inv(g)h))\in H^*$$
is an isomorphism. Therefore,  $\lambda$ is nondegenerate and $H$
is Frobenius.
\end{proof}
\section{Nichols-Zoeller for quasi-Hopf algebras}
In this section we prove the quasi-Hopf algebra version of
the Nichols-Zoeller theorem. The proof follows quite precisely
the course taken in the ordinary Hopf case, where we rely on the
account given in Montgomery's book \cite{Mon:HAAR}.
However, one has to take
special care whenever a Hopf module argument comes in, due to 
the fact that Hopf modules have a delicate definition. In particular,
the relevant Hopf module categories used in the original proof do 
not exist in the quasi-Hopf case; luckily, the actual objects in 
consideration can be located in Hopf module categories which do
exist.  

\begin{Thm}\nmlabel{Theorem}{auxthm}
  Let $K$ be a finite-dimensional quasi-Hopf algebra, and $W$ a 
  finitely generated right $K$-module. Suppose there exists a 
  finitely generated faithful right $K$-module $V$ such that 
  $W\o V\cong W^{\dim V}$ as right $K$-modules. Then $W$ is
  free over $K$.
\end{Thm}
\begin{proof}
  Since $K$ is Frobenius, we can apply
  \cite[Prop.3.1.2]{Mon:HAAR} which tells us that it is enough
  to show that some power of $W$ is free. Also, we can apply
  \cite[Lem.3.1.2]{Mon:HAAR} to conclude that there is $r>0$ 
  such that $W^r\cong F\oplus E$ with $F$ a free $K$-module 
  and $E$ a $K$-module that is not faithful, and also there is
  $s>0$ such that $V^s\cong F'\oplus E'$ with $F'$ free and
  $E'$ not faithful.
  Since $W^r\o V^s\cong (W\o V)^{rs}\cong W^{rs\dim V}$, we can 
  replace $W$ by $W^r$ and $V$ by $V^s$, and 
  assume that $W\cong F\oplus E$ and $V\cong F'\oplus E'$.

  So far, everything proceeded as in \cite{Mon:HAAR}. Now 
  we consider ${_\cdot^\cdot K_\cdot}\o {_\cdot V}\in\HMod K{}KK$
  by applying part \eqref{speziell} of \ref{remzeug} to $K^\bop$. By the version of
  \nmref{strukthm} for $K^\bop$, this is a free left $K$-module, hence 
  ${_\cdot K}\o {_\cdot V}\cong K^t$ as left $K$-modules, where
  $t:=\dim V$.
  Since $F$ is free, it follows that $F\o V\cong F^t$.
  Now
  $$F^t\oplus E^t\cong W^t\cong W\o V\cong (F\o V)\oplus (E\o V)\cong F^t\oplus (E\o V),$$
  and the Krull-Schmidt theorem implies
  $$E^t\cong E\o V\cong E\o(F'\oplus E')\cong(E\o F')\oplus (E\o E').$$
  Since $E$ is not faithful, neither is $E^t$, and hence neither is
  $E\o F'$. But $E\o F'$ is free, which finally implies that 
  $E=0$ and hence $W$ is free. To see that $E\o F'$ is free, it is in 
  turn sufficient to see that $E\o K$ is free, which follows from 
  \nmref{strukthm}, since we can consider
  ${E_\cdot}\o{_\cdot K^\cdot_\cdot}\in\HMod{}KKK$ by
  part \eqref{speziell} of \nmref{remzeug}.
\end{proof}
We arrive immediately at the quasi-Hopf version of the Nichols-Zoeller
theorem, for the more general situation advertised in the introduction.
\begin{Thm}\nmlabel{Theorem}{freeforalg}
  Let $H$ be a finite-dimensional quasi-Hopf algebra, and let
  $K\subset H$ be a subalgebra that admits a quasi-Hopf algebra
  structure. Then $H$ is a free (say, right) $K$-module.
\end{Thm}
\begin{proof}
  Consider the faithful right $K$-module $H$. By 
  \nmref{auxthm} it suffices to show that $H\o H\cong H^{\dim H}$ 
  as right $K$-module.
  But $H_\cdot\o{_\cdot H_\cdot^\cdot}\in\HMod{}HHH$ by
  part \eqref{cofree} of \nmref{remzeug}, 
  so that we can apply \nmref{strukthm} to conclude that
  $H\o H$ is indeed a free right $H$-module.
\end{proof}

To prepare for an application to semisimple quasi-Hopf algebras, 
modelled on \cite[Cor.3.2.3]{Mon:HAAR}, we will prove a result
of Panaite characterizing semisimple quasi-Hopf algebras in terms
of integrals; our proof will utilize Hopf modules.
\begin{Thm}[Panaite \cite{Pan:MTTQHA}]\nmlabel{Theorem}{Pan}
  A finite-dimensional quasi-Hopf algebra $H$ is semisimple if
  and only if it contains a left integral $t\in H$ with 
  $\epsilon(t)\neq 0$ (one might then as well assume that
  $\epsilon(t)=1$). 
  Here, an integral is an element satisfying
  $ht=\epsilon(h)t$ for all $h\in H$.
\end{Thm}
\begin{proof}
  If $H$ is semisimple, the epimorphism $\epsilon\colon H\rightarrow k$
  of left $H$-modules splits, which is the same as the 
  existence of 
  an integral $t$ with $\epsilon(t)=1$. 

  On the other hand, assume 
  that $\epsilon\colon H\rightarrow k$ splits as an $H$-module map.
  Then for any $V\in\LMod H$ the epimorphism 
  $H\o V\xrightarrow{\epsilon\o V} V$ splits. 
  ${_\cdot^\cdot H_\cdot}\o{_\cdot V}\in\HMod H{}HH$ is a Hopf module,
  hence a free left module by \nmref{strukthm}, which shows that
  $V$ is projective.
\end{proof}
\begin{Cor}
  Let $H$ be a finite-dimensional semisimple quasi-Hopf algebra. 
  Let $K\subset H$ be a subalgebra and (non-coassociative) 
  subcoalgebra 
  that has an associator making it a 
  quasi-Hopf algebra.

  Then $K$ is also semisimple.
\end{Cor}
\begin{proof}
  By \nmref{Pan}, $H$ contains a left integral $t$ with $\epsilon(t)\neq 0$.
  Since $H$ is a free left $K$-module, we can copy word-by-word
  the proof of \cite[Lem.2.2.2.2)]{Mon:HAAR} and conclude
  that $K$ also contains such an integral. By \nmref{Pan} again
  $K$ is semisimple.
\end{proof}

We now turn to the more general version of the
original Nichols-Zoeller 
theorem stating that every Hopf module in $\HMod{}H{}K$ is free
as a $K$-module. We shall show the same for such objects in 
$\HMod{}HKK$ that are finitely generated $K$-modules. In order
that the indicated Hopf modules be defined in the first place,
it will be necessary to assume that 
$K$ is a subquasibialgebra of $H$.
The weaker assumption in \nmref{freeforalg} would 
not be sufficient for this purpose.

\begin{Thm}
  Let $H$ be a finite-dimensional quasi-Hopf algebra and 
  $K\subset H$ a quasi-Hopf subalgebra. Then 
  every Hopf module $M\in\HMod{}HKK$ that is finitely generated
  as a $K$-module, is free as left $K$-module.
\end{Thm}
\begin{proof}
  By part \eqref{cofree} of \nmref{remzeug} 
  we can consider 
  ${_\cdot^\cdot H_\cdot}\o {_\cdot H_\cdot}\in\HMod H{}HH$, 
  and by \nmref{strukthm} it follows that 
  ${_\cdot^\cdot H_\cdot}\o {_\cdot H_\cdot}\cong {^\cdot_\cdot H_\cdot\o V_\cdot}$
  for some right $H$-module $V$. (In fact $V=H_{\operatorname{ad}}$, which is of
  no importance here.)

  Now we consider ${_\cdot^\cdot H_\cdot}\o{_\cdot H_\cdot}\in\HMod H{}KK$
  instead, where we still have the same isomorphism, of course.
  We conclude that, calculating with cotensor product within the 
  monoidal category $\BiMod K$ as indicated in 
  part \eqref{cotens} of  \nmref{remzeug}
  \begin{align*}
    {_\cdot M_\cdot}\o {_\cdot H_\cdot}
      &\cong{_\cdot M_\cdot^\cdot}\co H({_\cdot^\cdot H_\cdot}\o{_\cdot H_\cdot})\\
      &\cong M\co H({_\cdot^\cdot H_\cdot}\o V_\cdot)\\
      &\cong {_\cdot M_\cdot }\o V_\cdot
  \end{align*}
  Forgetting about the right module structure, we conclude
  $M\o H\cong M^{\dim V}$ as left $K$-modules, and apply 
  \nmref{auxthm} to conclude that $M$ is a free $K$-module.
\end{proof}

%\bibliographystyle{acm}
%\bibliography{promo}

\begin{thebibliography}{10}

\bibitem{Dri:QHA}
{\sc Drinfel'd, V.~G.}
\newblock Quasi-{H}opf algebras.
\newblock {\em Leningrad Math. J. 1\/} (1990), 1419--1457.

\bibitem{HauNil:DCPDQQG}
{\sc Hausser, F., and Nill, F.}
\newblock Diagonal crossed products by duals of quasi-quantum groups.
\newblock {\em Rev. Math. Phys. 11\/} (1999), 553--629.

\bibitem{HauNil:DQQG}
{\sc Hausser, F., and Nill, F.}
\newblock Doubles of quasi-quantum groups.
\newblock {\em Comm. Math. Phys. 199\/} (1999), 547--589.

\bibitem{HauNil:ITQHA}
{\sc Hausser, F., and Nill, F.}
\newblock Integral theory for quasi-{H}opf algebras.
\newblock {\em preprint\/} (math.QA/9904164).

\bibitem{Kap:B}
{\sc Kaplansky, I.}
\newblock Bialgebras.
\newblock University of Chicago Lecture Notes, 1975.

\bibitem{Kas:QG}
{\sc Kassel, C.}
\newblock {\em Quantum Groups}, vol.~155 of {\em GTM}.
\newblock Springer, 1995.

\bibitem{Mas:HAEC}
{\sc Masuoka, A.}
\newblock Hopf algebra extensions and cohomology.
\newblock {\em preprint\/} (2001).

\bibitem{Mas:CCBEAMPB}
{\sc Masuoka, A.}
\newblock Cohomology and coquasi-bialgebra extensions associated to a matched
  pair of bialgebras.
\newblock {\em Adv. Math.\/} (to appear).

\bibitem{Mon:HAAR}
{\sc Montgomery, S.}
\newblock {\em Hopf algebras and their actions on rings}, vol.~82 of {\em CBMS
  Regional Conference Series in Mathematics}.
\newblock AMS, Providence, Rhode Island, 1993.

\bibitem{NicZoe:HAFT}
{\sc Nichols, W.~D., and Zoeller, M.~B.}
\newblock A {H}opf algebra freeness theorem.
\newblock {\em Amer. J. Math. 111\/} (1989), 381--385.

\bibitem{ObeSch:UFGEI}
{\sc Oberst, U., and Schneider, H.-J.}
\newblock Untergruppen formeller {G}ruppen von endlichem {I}ndex.
\newblock {\em J. Algebra 31\/} (1974), 10--44.

\bibitem{Pan:MTTQHA}
{\sc Panaite, F.}
\newblock A {M}aschke-type theorem for quasi-{H}opf algebras.
\newblock In {\em Rings, Hopf algebras, and Brauer groups (Antwerp/Brussels,
  1996)}. Dekker, New York, 1998, pp.~201--207.

\bibitem{Par:NARMTI}
{\sc Pareigis, B.}
\newblock Non-additive ring and module theory {I}. {General} theory of monoids.
\newblock {\em Publ. Math. Debrecen 24\/} (1977), 189--204.

\bibitem{Par:NARMTII}
{\sc Pareigis, B.}
\newblock Non-additive ring and module theory {II}. {$\mathcal C$}-categories,
  {$\mathcal C$}-functors and {$\mathcal C$}-morphisms.
\newblock {\em Publ. Math. Debrecen 24\/} (1977), 351--361.

\bibitem{Sch:HMDQHA}
{\sc Schauenburg, P.}
\newblock Hopf modules and the double of a quasi-{H}opf algebra.
\newblock {\em Trans. Amer. Math. Soc.\/} (to appear).

\end{thebibliography}

\end{document}